\documentclass[11pt]{article}

\usepackage{amsfonts,amsmath,amssymb,graphics,bm}

\newcommand{\A}{\mathfrak{A}}

\newcommand{\R}{\mathbb{R}}

\newcommand{\E}{\textrm{E}}

\newcommand{\pr}{\mathbf{P}}

\newcommand{\sech}{\textrm{sech}}

\newcommand{\notthis}[1]{}

\newcommand{\inv}{^{-1}}

\begin{document}

\title{\bf Root-finding: from Newton to Halley and beyond}
\author{Richard J. Martin\footnote{Department of Mathematics, Imperial College London, South Kensington, London SW7 2AZ, UK} 
}
\maketitle

\begin{abstract}

We give a new improvement over Newton's method for root-finding, when the function in question is doubly differentiable. It generally exhibits faster and more reliable convergence. It can be also be thought of as a correction to Halley's method, as this can exhibit undesirable behaviour. 

\end{abstract}


The problem of finding the zero (more properly, a zero, though in context there may only be one) of a function is a common problem in all branches of applied mathematics.
To take a small subset, from the field of quantitative finance, here are some possible applications:
\begin{itemize}
\item
Finding the yield of a bond from its price;
\item
Finding the Black--Scholes implied volatility of an option from its price;
\item
In using the saddlepoint method \cite{Martin11b} for finding the tail probability of a random variable $Y$, i.e.\ $\pr(Y>y)$, given its cumulant generating function $K_Y$, we need to solve $K_Y'(s)=y$ for $s$;
\item
In generating random deviates using the inverse-cumulative method \cite[\S7.2]{NRC} we are to solve $F(y)=p\in(0,1)$.
\item
Computation of inverse special functions such as trigonometric or cumulative Normal is a ubiquitous problem. (That said, many of these problems are so standard that specialist optimised routines have already been written for the purpose.)
\end{itemize}
Without loss of generality we can rewrite these problems as $f(x)=0$ where $f$ is given, and often it will be the case that we can easily evaluate the first and second derivatives of $f$ without significant extra effort. The question is how best to employ this information.

\section{Root-finding methods}

\subsection{Newton and Extended Newton}

The Newton--Raphson formula \cite[\S9.4]{NRC} for finding a zero of a differentiable function $f$, also known simply as Newton's method, is
\begin{equation}
x \mapsto x - \frac{f(x)}{f'(x)}.
\label{eq:NR}
\end{equation}
It is exact, in the sense of converging in one step, when $f$ is linear, and its accuracy is second-order, i.e.\ in the vicinity of the root the error is proportional to the square of the previous error.

It is natural to ask how to incorporate information about higher derivatives. The following appears to be a new result, and it is exact whenever $f$ is the form $a\log(bx+c)$:
\begin{equation}
x \mapsto x - \frac{f(x)}{f'(x)} \, \mathcal{E}(q) ;
\qquad
q = \frac{f(x)f''(x)}{f'(x)^2},
\qquad
\mathcal{E}(t)\equiv (e^t-1)/t.
\label{eq:HNR1}
\end{equation}

\subsection{Comparison with Halley's method}

Halley's method, by contrast, is exact when $f$ is a M\"obius function (linear $\div$ linear).
Like ours it is third-order accurate when close to the root.
Its expression is similar:
\begin{equation}
x \mapsto x - \frac{f(x)}{f'(x)} \times \frac{1}{1-q/2} .
\end{equation}
However, it goes wrong when $q\ge2$.  Really, $x$ should always move in the opposite direction to $f(x)/f'(x)$, so that whatever multiplies $f(x)/f'(x)$ must be positive. It is easy to see where the fault lies: M\"obius functions do not map any domain of $\R$ onto $\R$ (onto = surjective).
The following example illustrates this: suppose $x=0$, $f(0)=1$, $f'(0)=1$. Then Halley approximates $f(x)\approx \frac{1+(c+1)x}{cx+1}$ where $c=-f''(0)/2$. All is well if $f''(0)<2$, as then the estimated root lies to the left of the origin, at $\frac{-1}{1+c}$; but if $f''(0)$ is any higher then it flips to being positive. The problem is that the range of $\frac{1+(c+1)x}{cx+1}$, for $x\le0$, does not include the origin in that case, so no sensible estimate of the root can be found.

Both $\mathcal{E}(q)$ and $(1-q/2)\inv$ have Maclaurin expansions $1+q/2+O(q^2)$, and indeed any higher-order correction of Newton's method must have this property.
This can be seen by writing $\xi$ for the root, stipulating that $\xi=x-[f'(x)/f(x)]C(q)$, where $C(q)$ is to be found, and then developing the RHS in a Taylor series around $x=\xi$.

Returning to (\ref{eq:HNR1}), it is worth noting that the presence of an exponential is not especially desirable: if $x$ is a long way away from the root, it can produce a correction step that is far too large or lead to an overflow error. A sort of `poor man's version' of (\ref{eq:HNR1}) is obtained using Pad\'e approximation, with different choices of rational approximation according as $q\gtrless0$:

\begin{equation}
x \mapsto x - \frac{f(x)}{f'(x)} \times \left\{
\begin{array}{ll} 1+(q/2)(1+q/3), & q\ge0 \\ 1/(1-q/2), & q\le 0 \end{array}
 \right.
\label{eq:HNR2}
\end{equation}
which preserves a link to Halley's method, avoids the fierceness of the exponential, and is quicker to compute. 
{\bf In all our numerical tests it is this version that we use.}

\subsection{The reasoning behind $a\log(bx+c)$}

We have shown why Halley's method can go wrong; our recipe will always put the next estimate on the correct side of the current one, in the sense that $x_{n+1}\gtrless x_n$ according as $f(x_n)f'(x_n)\lessgtr0$, because log continuously maps $\R^+$ onto $\R$. However, there are many functions that do that: we could use affine transformation of $x-x\inv$, for example. 

What is unique about the log function is that we are looking for a three-parameter family of functions, and it arises from Lie group theory in a natural way, as follows. Let $\A_x,\A_y$ be the (continuous, Lie) group of affine transformations in the $x$ and $y$ direction. For an arbitrary generating function---for example, $x-x\inv$ as suggested above---transformations under $\A_x\times \A_y$ will generate a four-parameter family of functions.
However, if $f$ is left invariant by some combination of actions (i.e.\ its stabiliser is non-trivial) then we will have only a three-parameter family---or even fewer, if we choose $f(x)=x$ as our generator. Writing $y=f(x)$ and applying the infinitesimal generator under which it is invariant, we obtain for some coefficients $c_j$:
\[
(c_1+xc_2)f'(x) = c_3+c_4f(x).
\]
This differential equation is easily solved. The only solution with the desired surjectivity property is the log function.

\section{Examples}

We give examples that compare with Newton's method. Some are taken from the current Wikipedia page on the subject. Empirical work suggests that the basin of attraction is larger, and sometimes much larger, than that of Newton's method.

\subsection*{Rapidity of convergence: $x^2-612=0$}

Starting from $x_0=10$, with Newton's method the iterates have errors ($x_n-\xi$) as follows\footnote{$x\E y$ means $x\cdot 10^y$; in context, `$<\!\epsilon$' means less than about $10^{-15}$. }:
\[
-1.47\E{1},\;
1.09\E{1}, \;
1.66\E{0},\;
5.20\E{-2},\;
5.45\E{-5},\;
6.01\E{-11}, \;
{<\!\epsilon}
\]
whereas with (\ref{eq:HNR2}) they are
\[
-1.47\E{1},\;
-3.51\E{0},\;
-2.20\E{-2},\;
-4.37\E{-9},\;
{<\!\epsilon}.
\]

\subsection*{Cycle: $\tanh x=0$}
If $|x_0| \gtrsim 1.088659$, this being the root of the equation $\tanh a = 2a \,\sech^2 a$, Newton's method is unstable. The iteration maps $a$ to $-a$ and vice versa.
Such behaviour is also seen with (\ref{eq:HNR2}) but the basin of attraction is $|x_0| \lesssim 2.410975$, over twice the size.

\subsection*{Turning-points: $x^3-2x+2=0$}
This is another example of a 2-cycle, as starting Newton's method at 0 or 1 will flip to the other point. Here, though, there is the additional problem of the turning-points at $\pm\sqrt{2\over3}$. 
With (\ref{eq:NR}) it is clear that at any turning-point the iteration is undefined. With (\ref{eq:HNR1},\ref{eq:HNR2}) the position is more complicated because it depends on the sign of $f''$ in relation to $f$.

With this particular example, if we start at the left turning-point $x=-\sqrt{2\over3}$, where $f>0$ but $f''<0$, we have $q=-\infty$ and so the iteration has an unstable stationary point there. To see why this is so, we approximate the iterative step as\footnote{The behaviour of $(\ref{eq:HNR1},\ref{eq:HNR2})$ is the same up to a factor of 2. The displayed equation relates to (\ref{eq:HNR2}).}
\[
\delta x \approx \frac{2f(x)}{f'(x)q} =  \frac{2f'(x)}{f''(x)}
\]
which is regular at the turning-point. So if we start to the left of $-\sqrt{2\over3}$, the iteration will reach the root, $\xi\approx -1.769292$.

On the other hand, if we start on the right, we can be led on a merry dance. Starting at $x=0.125$, it does find the root after 70 iterations, but in the meantime has gone as high as $2.4\times10^{10}$.
Yet starting at 0.0625 gives convergence to machine precision in 4 iterations. In essence the behaviour is governed by how close it gets to the other turning-point, $x=\sqrt{2\over3}$, where now $ff''>0$: the closer it gets, the further it is catapulted off. At a turning-point for which $ff''>0$, the iteration is undefined, because $q=+\infty$ then.

\subsection*{Newton fractal: $x^3-2x^2-11x+12=0$}
The roots are $-3,1,4$ and this serves as an illustration that the basins of convergence of Newton's method may be very small and give rise to intricate and interesting behaviour. For example \cite{Dence97}, if $x_0$ is moved from 2.3528363 up to 2.35287527 in small steps, different roots are found, as follows: 1, $-3$, 4, then $-3$ again, then 4 again. The behaviour of (\ref{eq:HNR2}) seems much less interesting. In this example, empirical investigation suggests that the boundaries are at around $-1.360920$ and 2.694254, and convergence is to the obvious root in each case.

\bibliographystyle{plain}
\bibliography{}

\end{document}